\documentclass[11pt]{amsart}
\usepackage{latexsym,amsmath,amsthm,amssymb}

\newcommand{\dR}{\ensuremath{\mathbb{R}}} 
\newcommand{\dZ}{\ensuremath{\mathbb{Z}}} 

\newcommand{\GI}{\mathbf{L}} 
\newcommand{\ind}{\mathrm{1}\hskip -3.2pt \mathrm{I}} 
\newcommand{\ent}{\mathbf{Ent}} 
\newcommand{\var}{\mathbf{Var}} 

\newtheorem{theorem}{Theorem}
\newtheorem{proposition}[theorem]{Proposition}
\newtheorem{lemma}[theorem]{Lemma}

\theoremstyle{definition}

\theoremstyle{remark} 
\newtheorem{remark}[theorem]{Remark}

\begin{document}   

\title{The Logarithmic Sobolev Constant of The Lamplighter}
\author[Abakumov et al.]
{Evgeny~Abakumov, Anne~Beaulieu, Fran\c{c}ois~Blanchard, Matthieu~Fradelizi, Natha\"el~Gozlan, Bernard~Host, Thiery~Jeantheau, 
Magdalena~Kobylanski, Guillaume~Lecu\'e, Miguel~Martinez, Mathieu~Meyer, Marie-H\'el\`ene~Mourgues, 
Fr\'ed\'eric~Portal, Francis~Ribaud, Cyril~Roberto, Pascal~Romon, Julien~Roth, 
Paul-Marie~Samson, Pierre~Vandekerkhove, Abdellah~Youssfi$^*$
}

\thanks{* The interested reader may also have a look to the following contribution by the authors: "Compter et mesurer. Le souci du nombre dans
l'\'evaluation de la production scientifique" [french], to appear in {\it La Gazette des Math\'ematiciens}, Soci\'et\'e Math\'ematique de France, 2010}

\date{\today}

\address{Universit\'e fran\c{c}aise - Laboratoire d'Analyse et de Math\'e\-matiques Appliqu\'ees (UMR CNRS 8050), 5 bd Descartes, 77454 Marne-la-Vall\'ee Cedex 2, France}
\email{email adresses follow the rule: name.surname@univ-mlv.fr}

\keywords{logarithmic-Sobolev inequalities, wreath-product, lamplighter graph, spectral gap}
\subjclass{60E15, 60J27, 60J45, 55C99 and 26D10}


\begin{abstract}
We give estimates on the logarithmic Sobolev constant of some finite lamplighter graphs in terms of the spectral gap of the 
underlying base. Also, we give examples of application.
\end{abstract}

\maketitle


Let $G =(V,E)$ be a finite graph with set of vertices $V$ and set of edges $E$. To each vertex of $G$ we associate
a lamp that can be on or off. Then, consider a lamplighter walking on $G$ and changing randomly the value of the lamps. The labelling configuration of the lamps together with the position of the lamplighter have a structure of graph known as \emph{lamplighter graph} or \emph{wreath product}. It is denoted by $G^\diamond=\{0,1\} \wr G$.

Mathematically, the vertices of $G^\diamond$ are all pairs $(\sigma,x) \in \{0,1\}^G \times G$ (the value of $\sigma(x) \in \{0,1\}$ gives the status of the lamp associated to site $x$: it is on if $\sigma(x)=1$, and off otherwise). Also, there is an edge between $(\sigma,x)$ and $(\sigma',x')$ in $G^\diamond$, if either $\sigma = \sigma'$ and $(x,x') \in E$ (which corresponds to saying that the lamplighter moved between $x$ and $x'$, on $G$, without changing the value of any lamp), or $\sigma(y)=\sigma'(y)$ for all $y \neq x$ and $x=x'$  (which corresponds to saying that the lamplighter did not move but might have changed the value of the lamp at site $x$). We shall say that $G$ is the \emph{base} of the lamplighter graph.

For example, the wreath product $(\mathbb{Z}_2^n)^\diamond=\{0,1\} \wr \mathbb{Z}_2^n$ of the two dimensional torus of side $n$ can be though of as the streets of a very regular city (as {\it e.g.} Buenos Aires) along which walks a lamplighter who switches the lamps on and off while passing by.

When $G$ is a group, the wreath product can be defined alternatively as the Cayley graph of some semi-direct product 
$G \ltimes \{0,1\}$ (see {\it e.g.} \cite{pittet-saloff-coste}). Note that our definition may slightly differ from others.
Also the wreath product can be defined more generally on an infinite base $G$, and also with lamps taking values in a more general graph $H$, leading to the wreath product $H \wr G$ (see {\it e.g.} \cite[Section 2]{pittet-saloff-coste}).
The geometry of the wreath product has many interesting properties, see for example 
\cite{Grigorchuk-zuk,vershik,erschler,gromov,austin-naor-peres} and references therein. On the other hand, many authors have studied random walks on the wreath product: the lamplighter performs a random walk on the base and changes randomly the value of the lamps, see
\cite{haggstrom,kaimanovich,pittet-saloff-coste,peres-revelle,tetali} and references therein.
 In particular, the wreath-product of the two-dimensional torus $(\mathbb{Z}_2^n)^\diamond$ is a very simple (and one of the rare) example of graph for which the relaxation time, the total variation mixing time, the entropy mixing time and the uniform mixing time were shown to be all of different orders of magnitude \cite{peres-revelle,tetali}.

In this note we shall estimate the logarithmic Sobolev constant of the lamplighter graph. In particular, we will prove that,
under mild assumption on the base, the logarithmic Sobolev constant of the wreath product can be expressed in terms of the spectral gap of the random walk on $G$ and on some other quantity involving only $G$. Our result can be extended in many ways (see Remark \ref{rem:extension}). However, for seek of clarity, we decided to focus on some specific random walks on $G^\diamond$ that will already provide interesting examples of application.

The logarithmic Sobolev (in short log-Sobolev) constant is a very powerfull tool in the study of high (and infinite) dimensional systems. It finds application in many areas of mathematics, see the books \cite{grossbook},
\cite{bakry}, \cite{Dav}, \cite{ane}, \cite{ledoux}, \cite{Wbook}, \cite{GZ99}, \cite{Ro99} for an introduction.
In particular it gives control on Gaussian concentration of Lipschitz functions and provides estimates on the convergence to equilibrium of stochastic processes. 

It is well known that the log-Sobolev constant gives a lower and an upper bound on the uniform mixing time, see 
\cite[Corollary 2.2.7]{saloff}. 
Since Peres-Revelle \cite{peres-revelle} and Ganapathy-Tetali \cite{tetali} obtained estimates on this quantity,
some bounds on the log-Sobolev constant can readily be derived from their papers. Our results however are sharper. Moreover, their approach is concerned with probabilistic tools (cover times, hitting times)
while  our analysis  is more on the potential theory side. As a consequence our results (even if not stated in the most general way) hold
in more general settings than the one considered in \cite{peres-revelle} and \cite{tetali}: in particular, we can consider a randomisation of the lamps that depends on the values of the lamps around, rising to
the Ising model at high temperature or even more general interacting particle systems for the lamps.
Finally, Peres-Revelle \cite{peres-revelle} and Ganapathy-Tetali \cite{tetali} have an extra hypothesis on the geometry of the base $G$ that is not necessary in our approach. 

We have to notice that in almost all the examples we will consider (except the hypercube for which the uniform mixing is not known), the log-Sobolev constant appears to be of the same order of magnitude of the uniform mixing time. This unfortunately makes it a weak tool in the analysis of the rate of convergence to equilibrium of the wreath-product.

Finally, we note that 
we are able to give an upper bound on the log-Sobolev
constant of the wreath product using a simpler functional inequality on the base (Poincar\'e inequality). We believe that our technique
may be one step in the direction of understanding how to deal with more sophisticated functional inequalities on the lamplighter graph
(as for instance Nash inequalities), using other, and hopefully simpler, functional inequalities on the base. In turn, this could maybe lead to an alternative proof of the very deep result of Erschler \cite{erschler}.

The paper is organized as follows. In the next section we give a precise definition of the log-Sobolev constant, of the dynamics and state our main results. Section 2 is dedicated to the proof of our results, while the last section deals with some examples.

\section{Setting and results}
In order to introduce the Markov generator of the dynamics, we need first to give some notations.
Given a configuration $\sigma \in \{0,1\}^G$ ({\it i.e.} a labelling of the lamps), and a site $x \in G$, we state
$\sigma^x$ for the configuration flipped at site $x$:
$$
\sigma^x (y) = \left\{ 
\begin{array}{ll}
\sigma(y) & \mbox{if } y \neq x \\
1 - \sigma(x) &  \mbox{if } y = x .
\end{array}
\right.
$$
We will denote by $p(x,y)$ the transition probabilities for the lamplighter to move from $x$ to $y$, and by $c_x(\sigma)$
the transition rate from the configuration $\sigma$ to $\sigma^x$.
Also, given two sites $x$, $y$ of $G$, the notation $x \sim y$ means that $x$ and $y$ are nearest neighbors in $G$, {\it i.e.} $(x,y) \in E$.

Then, the Markov process of the wreath product we are interested in is given by the following infinitesimal generator,
acting on any function $f : G^\diamond \to \mathbb{R}$,
$$
\GI f (\sigma,x) = \frac{1}{2} \GI_b f (\sigma,x) + \frac{1}{2} \GI_\ell f(\sigma,x)
$$
with
$$
\GI_b f (\sigma,x) = \sum_{y \sim x} p(x,y) \left( f(\sigma,y) - f(\sigma,x) \right)
$$
and
$$
\GI_b f (\sigma,x) = c_x(\sigma) \left(f(\sigma^x,x) - f(\sigma,x) \right) .
$$
The operator $\GI_b$ corresponds to the dynamics of the lamplighter on the base, while $\GI_\ell$
corresponds to the dynamics of the lamps.

We assume that the transition probabilities/rates are reversible with respect to some probability measure $\nu$ on $G$ for the lamplighter,
and with respect to some probability measure $\mu$ on $\{0,1\}^G$ for the lamps.
In other words, for any $x$, $y$ in $G$, and any $\sigma \in \{0,1\}^G$, the following detailed balance conditions hold
$$
\nu(x) p(x,y) = \nu(y) p(y,x)  \quad \mbox{and} \quad \mu(\sigma)c_x(\sigma) = \mu(\sigma^x) c_x(\sigma^x) .
$$
Now define $\pi = \mu \times \nu$. The above equalities guarantee that the generator $\GI$ is symmetric in 
$\mathbb{L}_2(\pi)$ or equivalently, that the Markov process is reversible with respect to $\pi$.

A probabilistic interpretation of the generator $\GI$ is the following. Wait a mean $1$ exponential time and then launch a coin toss. If the coin is head, then the lamplighter moves at random to one of its neighbors according to the transitions given by $p$, and the lamps do not change. If instead the coin is tail then the lamplighter does not move but randomizes the lamp of his/her position
(say $x$) according to the transition given by $c_x (\sigma)$. Then the procedure starts again.

While the reversible measure $\pi$ is product, which usually makes life simpler, the dynamics is not, since a lamp can be randomized only when the lamplighter is sitted on the corresponding site. The main difficulty in our analysis will come from this non-product character.

Finally we introduce the Dirichlet form $\mathcal{E}_\pi$ of the wreath product
$G^\diamond$, that is the non-negative bilinear form $\mathcal{E}_\pi(f,g) = - \pi(f \cdot \GI g)$ acting on functions
$f,g : G^\diamond \to \mathbb{R}$ (the dot sign corresponds to the scalar product in $\mathbb{L}_2(\pi)$). Namely, using
the detailed balance condition above, one can write
\begin{align*}
4 \mathcal{E}_\pi (f,g) & = 
\!\!\!\! \sum_{\sigma \in \{0,1\}^G} \sum_{x \in
G} \pi(\sigma,x) \sum_{y \sim x} p(x,y) \left[f(\sigma,y) - f(\sigma,x)\right]\left[g(\sigma,y) - g(\sigma,x)\right] \\
& + \!\!\sum_{\sigma \in \{0,1\}^G} \sum_{x \in
G} \pi(\sigma,x) c_x(\sigma) \left[f(\sigma^x ,x) - f(\sigma,x)\right]\left[g(\sigma^x ,x) - g(\sigma,x)\right]  .
\end{align*}
For simplicity we set $\mathcal{E}_\pi (f) = \mathcal{E}_\pi (f,f)$.

The main quantity of interest for us is the logarithmic Sobolev constant $C_{LS}(\pi)$
which is defined as the best constant $C$ satisfying for any function
$f:G^\diamond \rightarrow \dR$,
\begin{equation} \label{eq:ls}
\ent_\pi (f^2) \leq C \mathcal{E}_\pi (f)
\end{equation}
where $\ent_\pi (f^2) := \pi(f^2 \log(f^2/ \pi(f^2)))$ is the entropy
of $f^2$ and $\pi(g)$ is a short hand notation for the mean value of 
$g$ under $\pi$.

In order to state our results, we also need to define the spectral gap of $(G,\nu)$. 
Let
\begin{equation} \label{eq:enu}
\mathcal{E}_\nu (f,g)  :=
\frac 12  \sum_{x \in G} \sum_{y \sim x} \nu(x) p(x,y) \left[f(y) - f(x)\right]\left[g(y) - g(x)\right]
\end{equation}
be the Dirichlet form of $G$, acting on functions $f,g : G \to \mathbb{R}$. Here also we set 
$\mathcal{E}_\nu (f)=\mathcal{E}_\nu (f,f)$.
Then, the spectral gap constant ${\rm gap}(\nu)$ of $(G,\nu)$ is the best constant $C$ such that for any $f$, 
$$
C \var_\nu(f) \leq \mathcal{E}_\nu (f)
$$
where $\var_\nu(f)=\nu(f^2) - \nu(f)^2$ and again $\nu(g)=\sum_{x \in G} \nu(x) g(x)$ is a short hand notation for the mean value of $g$ under $\nu$.

We are in position to give our first main result. Set $\nu^*:=\min_{x \in G} \nu(x)$.

\begin{theorem} \label{th:main-intro}
Assume that $\mu$ is the unifrom measure on $\{0,1\}^G$ and that for any $x$ and any $\sigma$, $c_x(\sigma) \leq 1/2$.
Then,
$$
C_{LS}(\pi) \leq \frac{6}{\nu^* {\rm gap}(\nu)}.
$$
\end{theorem}

The two assumptions of the theorem were also considered in \cite{peres-revelle,tetali}.
Moreover, in \cite{peres-revelle}, it is assumed that $c_x(\sigma)=\frac{1}{2}p(x,x)$. This justifies the 
choice of the value $1/2$.

When $\nu \equiv \frac{1}{|G|}$, where $|G|$ denotes the cardinality of $G$, 
we have also a lower bound on the log-Sobolev constant.
We will give the result for a particular class of sets we introduce now.
For any subset $B$ of $G$, we define its closure
by $\bar{B}:=\{x \in G : \exists y \in B \mbox{ such that } x \sim y \}$.
We will say that $G$ satisfies Hypothesis $(H)$ with parameter $\varepsilon$
if there exists $\varepsilon >0$
such that
$$
|\bar B| \geq \frac 12 |G| \Rightarrow |B| \geq \varepsilon |G| .
$$
This assumption guarantees that the points of $G$ do not have too many neighbors.

\begin{proposition} \label{prop:low}
Assume that $G$ satisfies Hypothesis $(H)$ with parameter $\varepsilon$ and
that $\mu$ and $\nu$ are the uniform measures on $\{0,1\}^G$ and $G$ respectively (in particular $\nu^* \equiv 1/|G|$).
Then,
$$
C_{LS}(\pi) \geq \varepsilon  \frac{|G|}{{\rm gap}(\nu)} .
$$
\end{proposition}

\begin{remark}
Under the assumptions of Theorem \ref{prop:low}, and assuming for instance that $c_x(\sigma) = 1/2$ for all $x$ and all $\sigma$,
the two bounds of Theorem \ref{th:main-intro} and Proposition \ref{prop:low} match. Namely
$$
 \varepsilon \frac{|G|}{{\rm gap}(\nu)}  
\leq C_{LS}(\pi)
\leq
6 \frac{|G|}{{\rm gap}(\nu)}  .
$$
\end{remark}

\section{Proofs}
This section is dedicated to the proofs of our main theorems. We shall start with Theorem \ref{th:main-intro}
that we state in a more general form.

\subsection{Proof of Theorem \ref{th:main-intro}}
In order to state Theorem \ref{th:main-intro} in a more general form, we need to introduce the logarithmic Sobolev constants of $(G,\nu)$ and of $(\{0,1\}^G,\mu)$. Let
$$
\mathcal{E}_\mu (f) :=
\frac 12 \sum_{\sigma \in \{0,1\}^G} \sum_{x \in
G} \mu(\sigma) c_x(\sigma) (f(\sigma^x) - f(\sigma))^2
$$
be the Dirichlet form of $(\{0,1\}^G,\mu)$ (acting on functions $f:\{0,1\}^G \rightarrow \dR$) and recall the definition of $\mathcal{E}_\nu$ given in \eqref{eq:enu}.

Then, we denote by $C_{LS}(\mu)$ and $C_{LS}(\nu)$
the logarithmic Sobolev constants of $(\{0,1\}^G,\mu)$ and $(G,\nu)$ respectively, namely
the best constant $C$ and $C'$ such that for any
$f:\{0,1\}^G \rightarrow \dR$,
$$
\ent_\mu (f^2) \leq C \mathcal{E}_\mu (f) ,
$$
respectively for any $f : G \rightarrow \dR$,
$$
\ent_\nu (f^2) \leq C' \mathcal{E}_\nu (f) .
$$

We shall prove the following.

\begin{theorem} \label{th:main}
Let $\nu^*:=\min_{x \in G} \nu(x)$ and assume that $c_x(\sigma) \leq a$ for some constant
$a$ independent on $x\in G$ and $\sigma \in \{0,1\}^G$.
It holds,
$$
C_{LS}(\pi) \leq \frac{1}{2\nu^*} \max \left( \frac{1+12a} {{\rm gap}(\nu)} , 6C_{LS}(\mu) \right) .
$$
\end{theorem}

We postpone the proof of Theorem \ref{th:main} in order to derive the result of Theorem \ref{th:main-intro}.

\begin{proof}[Proof of Theorem \ref{th:main-intro}]
Note first that, under the assumptions of Theorem \ref{th:main-intro}, $a \leq 1/2$. Moreover,
since $\mu$ is the uniform measure
on $\{0,1\}^G$, it is well known that $C_{LS}(\mu) = 2$ (see \cite[Chapter 1]{ane}, \cite{bonami}). Hence,
Theorem \ref{th:main} implies that
$$
C_{LS}(\pi) \leq \frac{1}{2\nu^*}
\max \left(
\frac{7}{{\rm gap}(\nu)} , 12 \right).
$$
The expected result follows, since ${\rm gap}(\nu) \leq 1$.
\end{proof}

Now we turn to the proof of Theorem \ref{th:main}.

\begin{proof}[Proof of Theorem \ref{th:main}]
Let $\xi \in G$ be the random variable giving the position of the lamplighter on $G$.
By the standard conditioning formula of the entropy, for any $f:G^\diamond \rightarrow \dR$
\begin{equation} \label{eq:start}
\ent_\pi(f^2) =
\ent_\pi \left( \pi(f^2 | \xi) \right)
+
\pi \left( \ent_{\pi(\cdot|\xi)}(f^2) \right) ,
\end{equation}
where as usual $\pi(\cdot|\xi)$ is the expected measure, given $\xi$.

First we deal with the first term in the right hand side of the latter. Since $\pi(f^2|\xi)$ is a function on $G$,
by the logarithmic Sobolev inequality of $(G,\nu)$, we have
\begin{align*}
&\ent_\pi \left( \pi(f^2 | \xi) \right)
 =
\ent_\nu \left( \pi(f^2 | \xi) \right) \\
&\quad \leq
\frac 12 C_{LS}(G,\nu) \sum_{x \in G} \sum_{y \sim x} \nu(x) p(x,y)
\left( \sqrt{\pi(f^2|\xi=y)}  - \sqrt{\pi(f^2|\xi=x)} \right)^2 .
\end{align*}
Write
$$
\sqrt{\pi(f^2|\xi=y)}  - \sqrt{\pi(f^2|\xi=x)}
=
\frac{\pi(f^2|\xi=y)  - \pi(f^2|\xi=x)}{\sqrt{\pi(f^2|\xi=y)} + \sqrt{\pi(f^2|\xi=x)}} 
$$
and notice that
\begin{eqnarray*}
\pi(f^2|\xi=y)  - \pi(f^2|\xi=x)
 =
\mu \left( \left[f(\cdot,y) - f(\cdot,x)\right]\left[f(\cdot,y) + f(\cdot,x)\right]\right) ,
\end{eqnarray*}
where by definition, $\mu \left(g(\cdot,x) \right):=
\sum_\sigma \mu(\sigma) g(\sigma,x)$ for any function $g$, any $x$.
Thus, by Cauchy-Schwarz inequality
\begin{align*}
& \left[ \sqrt{\pi(f^2|\xi=y)} - \sqrt{\pi(f^2|\xi=x)} \right]^2 \!
=
\frac{ \mu \left(\left[f(\cdot,y) - f(\cdot,x)\right]\left[f(\cdot,y) + f(\cdot,x)\right]
\right)^2 }{\left(\sqrt{\pi(f^2|\xi=y)} + \sqrt{\pi(f^2|\xi=x)} \right)^2} \\
& \qquad \leq
\frac{\mu\left( \left[f(\cdot,y) - f(\cdot,x) \right]^2 \right)
\mu\left( \left[f(\cdot,y) + f(\cdot,x)\right]^2\right)}
{\left(\sqrt{\pi(f^2|\xi=y)} + \sqrt{\pi(f^2|\xi=x)} \right)^2} \\
&  \qquad \leq
\sum_\sigma \mu(\sigma) (f(\sigma,y) - f(\sigma,x))^2
\end{align*}
since
\begin{align*}
& \frac{\mu\left( \left[f(\cdot,y) + f(\cdot,x)\right]^2\right)}
{\left(\sqrt{\pi(f^2|\xi=y)} + \sqrt{\pi(f^2|\xi=x)} \right)^2} \\
& \quad \qquad =
\frac{ \mu(f^2(\cdot,x)) + \mu(f^2(\cdot,y)) + 2 \mu(f(\cdot,x))\mu(f(\cdot,y))}
{\mu(f^2(\cdot,x)) + \mu(f^2(\cdot,y)) + 2 \sqrt{\mu(f(\cdot,x))\mu(f(\cdot,y))}}
\leq 1 .
\end{align*}
Summing up we arrive at
\begin{equation}\label{eq:step1}
\ent_\pi \left( \pi(f^2 | \xi) \right)
\leq
\frac 12 C_{LS}(G,\nu) \sum_\sigma \sum_{x \in G} \pi(\sigma,x)
\sum_{y \sim x}  p(x,y)
(f(\sigma,y) - f(\sigma,x))^2 .
\end{equation}
This quantity is the contribution of the random walk part performed by the lamplighter in the whole Dirichlet form.

For the second term in the right hand side of \eqref{eq:start} we use first the logarithmic Sobolev inequality
of $(\{0,1\}^G,\mu)$. Namely, we have
\begin{align*}
\pi \left( \ent_{\pi(\cdot|\xi)} (f^2) \right)
& =
\sum_{y \in G} \nu (y) \ent_\mu (f^2(\cdot,y)) \\
& \leq
\frac 12 C_{LS}(\mu)
\sum_{y \in G} \nu (y) \sum_\sigma \sum_{x \in G} \mu(\sigma) c_x(\sigma)
(f(\sigma^x,y) - f(\sigma,y))^2 \!.
\end{align*}
We decompose each term in the following telescopic sum
\begin{align*}
f(\sigma^x,y) - f(\sigma,y)
& =[f(\sigma^x,y)-f(\sigma^x,x)]
+
[f(\sigma^x,x)- f(\sigma,x)] \\
& \quad  +
[f(\sigma,x)- f(\sigma,y)] .
\end{align*}
Then, by Cauchy-Schwarz inequality, we get that
$$
\pi \left( \ent_{\pi(\cdot|\xi)} (f^2) \right)
\leq \frac 32 C_{LS}(\mu) \left( I+II+III \right)
$$
where
$$
I:=\sum_{y \in G} \nu (y) \sum_\sigma \sum_{x \in G} \mu(\sigma) c_x(\sigma)
(f(\sigma^x,y) - f(\sigma^x,x))^2 ,
$$
$$
II  := \sum_{y \in G} \nu (y) \sum_\sigma \sum_{x \in G} \mu(\sigma) c_x(\sigma)
(f(\sigma^x,x) - f(\sigma,x))^2
$$
and
$$
III:=\sum_{y \in G} \nu (y) \sum_\sigma \sum_{x \in G} \mu(\sigma) c_x(\sigma)
(f(\sigma,x) - f(\sigma,y))^2 .
$$
Using the detailed balance condition on $\mu$ and the change of variable
$\sigma \rightsquigarrow \sigma^x$ allows us to conclude that $I=III$. Moreover,
by Poincar\'e inequality for $(G,\nu)$,
\begin{align*}
III 
& \leq  
\frac{a}{\nu^*}
\sum_\sigma \mu(\sigma) \sum_{x \in G}  \sum_{y \in G}  \nu(x) \nu(y)
(f(\sigma,x) - f(\sigma,y))^2 \\
& = 
\frac{2a}{\nu^*}
\sum_\sigma \mu(\sigma) \var_\nu \left( f(\sigma,\cdot) \right) \\
& \leq 
\frac{a}{\nu^*{\rm gap}(\nu)}
\sum_\sigma \mu(\sigma)\sum_{x \in G}\sum_{y \sim x} \nu(x)p(x,y)
(f(\sigma,x) - f(\sigma,y))^2 .
\end{align*}
On the other hand,
\begin{align*}
II 
& =  
\sum_\sigma \sum_{x \in G} \mu(\sigma) c_x(\sigma)
(f(\sigma^x,x) - f(\sigma,x))^2 \\
& \leq 
\frac{1}{\nu^*}
\sum_\sigma \sum_{x \in G} \nu(x) \mu(\sigma) c_x(\sigma)
(f(\sigma^x,x) - f(\sigma,x))^2 .
\end{align*}

All the previous computations lead to
\begin{eqnarray*}
&& \pi \left( \ent_{\pi(\cdot|\xi)} (f^2) \right)
\leq
\frac  32 C_{LS}(\mu) \Big( \\
&& \quad \frac{2a}{\nu^* {\rm gap}(\nu)}
\sum_\sigma \mu(\sigma)\sum_{x \in G}\sum_{y \sim x} \nu(x)p(x,y)
(f(\sigma,x) - f(\sigma,y))^2 \\
&& \quad +
\frac{1}{\nu^*} \sum_\sigma \sum_{x \in G} \nu(x) \mu(\sigma) c_x(\sigma)
(f(\sigma^x,x) - f(\sigma,x))^2 \Big) .
\end{eqnarray*}
It follows from 
the very definition of $C_{LS}(\pi)$ that
$$
C_{LS}(\pi) \leq 
\max \left(C_{LS}(\nu) + \frac{6a}{\nu^* {\rm gap}(\nu)} , \frac{3}{\nu^*} C_{LS}(\mu) \right) .
$$
In order to conclude, we use a result by Diaconis and Saloff-Coste that compares the log-Sobolev constant to the spectral gap. Namely, they proved that $C_{LS}(\nu) \leq {\rm gap}^{-1}(\nu) \log \frac{1}{\nu^*}$  (see
\cite[Corollary 2.2.10]{saloff}, \cite{diaconis-saloff-coste}). Thus, we end up with
$$
C_{LS}(\pi) \leq 
\frac{1}{\nu^*} \max \left( \frac{\nu^* \log \frac{1}{\nu^*} + 6a} {{\rm gap}(\nu)} , 3C_{LS}(\mu) \right) .
$$
This ends the proof, since $|x \log x| \leq 1/2$ for $x \in (0,1)$. 
\end{proof}

\begin{remark}\label{rem:extension}
Our approach in Theorem \ref{th:main} is based on a very simple analysis of the position of the lamplighter on the base.
The conditioning property of the entropy allows to separate the effect of the random walk of the lamplighter himself/herself, and
the effect of the lamps. Only the product structure of $\pi=\mu \times \nu$ is used. This allows to consider more sophisticated
measure $\mu$ than the uniform measure on $\{0,1\}^G$. For example, the lamps can be randomized according to an interacting equilibrium measure such as the Ising model. See section \ref{sec:torus} for an illustration.

Also, one could generalized the proof by considering more values for the lamps, and a wreath product $G^\diamond = H \wr G$.
In this case, the cardinality of $H$ would certainly enters into the game.

Finally, one could consider a different kind of lamplighter. Assume for example that the number of lights that are on is fixed once for all.
Then, when the lamplighter turns on/off a lamp, he/she also has to turn off/on an other lamp. If this  one is of its nearest neighbors,
we end up with the Kawasaki dynamics on the lamps. Here the results would be very different in the application since
the Kawasaki dynamics is known to be very slow (the log-Sobolev constant is of order $n^2$ in the example of the
two (and actually $d$) dimensional torus of side $n$ (see \cite{lee-yau,yau,cancrini}).
\end{remark}

\subsection{Proof of Proposition \ref{prop:low}}
In order to prove Proposition \ref{prop:low}, we need first to introduce the notion of spectral profile.
Given a non-empty subset $S \subset G$, we set
$$
\lambda(S)= \inf_{f \in c_0(S)} \frac{\mathcal{E}_\nu(f)}{\var_\nu(f)}
$$
where $c_0(S)=\left\{ f \geq 0, \mathrm{supp}(f) \subset S, f \neq \mathrm{constant} \right\}$.
Then the spectral profile is defined by
\begin{equation} \label{eq:goel}
\Lambda(r) = \inf_{S:\nu(S) \leq r} \lambda(S) .
\end{equation}
This notion is related to Faber-Krahn inequalities introduced by Grigor'yan and developed with Coulhon and others
to estimate the rate of convergence of the heat kernel on manifolds, see \cite{grigoryan,coulhon}. The main interest for us is the following result by Goel, Montenegro and Tetali \cite[Lemma 2.2]{goel} relating the spectral profile to the spectral gap.

\begin{lemma}[\cite{goel}] \label{prop:goel}
It holds
$$
\frac{1}{\mathrm{gap}(\nu)} \leq \Lambda(1/2) \leq \frac{2}{\mathrm{gap}(\nu)} .
$$
\end{lemma}

We are now in position to prove Proposition \ref{prop:low}.

\begin{proof}[Proof of Proposition \ref{prop:low}]
Consider a function $g : G \rightarrow \dR$ with $\nu(g)=0$ and a set $A \subset G$
with $|\bar A| \geq |G|/2$ and the support of $g$ satisfying
 ${\rm supp}(g) \subset G \setminus {\bar A}$.
Denote by $\ind_{\sigma_A \equiv 1}$ the indicator function of the
configuration identically equal to $1$ on $A$.
Finally, let $f(\sigma,x)=g(x) \ind_{\sigma_A \equiv 1}(\sigma)$.

Note that $\pi(f^2)= \nu(g^2) \mu(\{\sigma_A \equiv 1\})$. Thus,
\begin{eqnarray*}
\ent_\pi(f^2) & = &
\sum_{\sigma : \sigma_A \equiv 1} \mu(\sigma) \sum_{x\in G} \nu(x) g^2(x)
\log \frac{g^2(x)}{\nu(g^2) \mu(\{\sigma_A \equiv 1\})} \\
& = &
\mu(\{\sigma_A \equiv 1\}) \left( \ent_\nu(g^2) + \nu(g^2)
\log \frac{1}{\mu(\{\sigma_A \equiv 1\})} \right) \\
& \geq &
\nu(g^2)\mu(\{\sigma_A \equiv 1\})\log \frac{1}{\mu(\{\sigma_A \equiv 1\})} .
\end{eqnarray*}
In the last line we used the well known fact that the entropy is non negative.

On the other hand, by construction, we have
\begin{eqnarray*}
 \mathcal{E}_\pi (f) & = &
\frac 14 \sum_{\sigma : \sigma_A \equiv 1} \mu(\sigma)
\sum_{x \in G} \sum_{y \sim x} \nu(x) p(x,y)(g(y)-g(x))^2 \\
& = &
\frac 12 \mu(\{\sigma_A \equiv 1\})  \mathcal{E}_\nu (g) .
\end{eqnarray*}

It follows by definition of the logarithmic Sobolev constant that
$$
C_{LS}(\pi) \geq 2 \frac{\nu(g^2)}{\mathcal{E}_\nu (g)}
\log \frac{1}{\mu(\{\sigma_A \equiv 1\})} .
$$
Since $\mu$ is uniform, by Hypothesis $(H)$
we have
$$
\mu(\{\sigma_A \equiv 1\})
= 
\frac{1}{2^{|A|}}
\leq 
\frac{1}{2^{\varepsilon |G|}} .
$$
And we are left with the bound
$$
C_{LS}(\pi) \geq 2 \varepsilon \log \frac{1}{2} |G|
\sup_{g: |{\rm supp}(g)| \leq \frac{|G|}{2}} \frac{\nu(g^2)}{\mathcal{E}_\nu (g)} 
\geq 
\varepsilon |G| \Lambda(1/2) 
$$
where $\Lambda(1/2)$ is the spectral profile introduced in \eqref{eq:goel}.
The expected result follows from Lemma \ref{prop:goel}.
\end{proof}

\begin{remark} \label{rem:1}
Note that, using test functions of the form $f : (\sigma,x) \mapsto g(\sigma)$, where $g : \{0,1\}^G \to \mathbb{R}$, respectively 
$f : (\sigma,x) \mapsto h(x)$, where $h : G \to \mathbb{R}$, we immediately get the following trivial bound
$$
C_{LS}(\pi) \geq 2 \max \left( C_{LS}(\mu), C_{LS}(\nu) \right) .
$$
\end{remark}

%
%
%


\section{Examples}

In this section we apply our results on some examples: the torus, the regular tree, the complete graph and finally the hypercube.
The first two examples will satisfy Hypothesis $(H)$, while the last two will not. In those cases, in order to estimate the
log-Sobolev contant, one has to consider accurate test functions in \eqref{eq:ls}. All our results are sharp, in the sens that
the log-Sobolev constant is shown to be of the right order of magnitude.

\subsection{The torus} \label{sec:torus}
Consider $G=\dZ_n^d$, the $d$ dimensional torus of side $n$. Assume that the lamplighter
performs a simple random walk on $G$. This implies that $\nu$ is the uniform measure: $\nu \equiv 1/n^d$.
Then, since each point of the torus has $2d$ neighbors, for any set $A$,
$|\bar A| \leq (1+2d) |A|$. In particular, Hypothesis $(H)$ holds
with $\varepsilon = \frac{1}{2(2d+1)}$. If $\mu$ is the uniform measure on $\{0,1\}^{\dZ_n^d}$, 
{\it i.e.} $\mu \equiv \frac{1}{2^{n^d}}$, since it is well known that 
${\rm gap}^{-1}(\nu) = O(n^2)$, it follows form our estimates that
$$
C_{LS}(\pi) = O(n^{d+2}) .
$$

The previous estimates can be easily adapted to the celebrated "dog" graph considered and analysed in \cite{diaconis-saloff-coste}
(we omit details).

Now consider the Ising measure $\mu$ at high temperature. It is know that $C_{LS}(\mu)$ is of order $1$ (\cite{martinelli}).
Hence, Theorem \ref{th:main} applies and leads to $C_{LS}(\pi) \leq O( n^{d+2})$. The strategy of the proof of
Proposition \ref{prop:low} also applies, with some minor modifications (we omit details), leading to a lower bound of the same type, and finally to $C_{LS}(\pi) = O(n^{d+2})$.

Finally consider the two point torus $G=\{0,1\}$ and $G^\diamond = \{0,1\} \wr G$, composed of $8$ points. Assume that $\nu(0)=\nu(1)=1/2$ and that $p(0,1)=p(1,0)=1/2$. Also, assume that $\mu$ is the product of Bernoulli$-p$ measures
and that, for any $x \in G$, $c_x(\sigma) = 1-p$ if $\sigma(x)=p$ and $c_x(\sigma) = p$ if $\sigma(x)=1-p$. Set $q=1-p$.
In \cite{diaconis-saloff-coste}, Diaconis and Saloff-Coste give the exact value of the log-Sobolev constant of $\mu$: $C_{LS}(\mu)=
\frac{\log(p) - \log(q)}{p-q}$, where $C_{LS}(\mu)=2$ (the limit) when $p=q=1/2$. See \cite{saloff,ane} for a simple proof of this fact, due to Bobkov. The log-Sobolev constant of the lamplighter graph is not easy to compute exactly, due to the $8$ points. Anyway, Theorem \ref{th:main} applies. Since $\nu^*=1/2$ and $\mathrm{gap}(\nu)=1$ (see 
\cite{ane}), we end up with
$$
C_{LS}(\pi) \leq 12 \frac{\log(p) - \log(q)}{p-q} .
$$
On the other hand, Remark \ref{rem:1} asserts that $C_{LS}(\pi) \geq C_{LS}(\mu)= \frac{\log(p) - \log(q)}{p-q}$.
Hence
$$
\frac{\log(p) - \log(q)}{p-q} \leq C_{LS}(\pi) \leq 12 \frac{\log(p) - \log(q)}{p-q}.
$$
In turn, our estimates catch (part of) the non-trivial behavior of $C_{LS}(\pi)$ when, say, $p \to 0$.

\subsection{Regular tree}
Consider $G=T$ the finite $b$-ary tree.
Assume that the lamplighter perform a simple
random walk on $T$ (which implies that $\nu \equiv 1/|T|$).
Assume that $\mu$ is the uniform measure on $\{0,1\}^T$. Note that each point
of $G$ has $b+1$ neighbors in such a way that hypothesis $(H)$ holds
with $\varepsilon = \frac{1}{2(b+2)}$. Since ${\rm gap}^{-1}(\nu) = O(|T|)$, the previous
results leads to
$$
C_{LS}(\pi) = O(|T|^2) .
$$
The above result can be obviously generalized to any graph with degree uniformly bounded from above and below.

\subsection{The complete graph}
Consider $G = K_n$ the complete graph of cardinal $n$: each point has $n$ neighbors.
Consider the simple random walk on $K_n$ reversible with respect to the uniform measure $\nu$ ($\nu \equiv 1/n$) with transition matrix $p(x,y)=1/n$. 
It is known that ${\rm gap}^{-1}(\nu)=1$.
Consider the uniform measure $\mu$ on $\{0,1\}^{K_n}$: $\mu \equiv 1/2^n$. Also, assume that $c_x(\sigma)=1/2$ for all $x$ and $\sigma$.
Theorem \ref{th:main-intro} leads to
$$
C_{LS}(\pi) \leq 6n .
$$
On the other hand, Hypothesis $(H)$ is not satisfied. Hence, in order to get a lower on the log-Sobolev constant, we have to find 
an accurate test function in \eqref{eq:ls}. Set $f(\sigma,x)=\ind_{\sigma \equiv 1}$, the indicator function that all the lamps are all
on. Then, 
\begin{align*}
\mathcal{E}_\pi(f) 
&= 
\frac{1}{4} \sum_\sigma \sum_{x \in K_n} \pi(\sigma,x) c_x(\sigma) \left[ \ind_{\sigma \equiv 1} (\sigma^x)- \ind_{\sigma \equiv 1}(\sigma) \right]^2 \\
&=
\frac{1}{8n} \sum_\sigma \sum_{x \in K_n} \mu(\sigma) \left[ \ind_{\sigma \equiv 1} (\sigma^x)- \ind_{\sigma \equiv 1}(\sigma) \right]^2 \\
& =
\frac{1}{4} \mu(\sigma \equiv 1) .
\end{align*}
On the other hand, since $\ent_\pi(f^2)=\ent_\mu(f^2) = \mu(\sigma \equiv 1) \log \frac{1}{\mu(\sigma \equiv 1)}$, 
the log-Sobolev inequality \eqref{eq:ls} reads
$$
\mu(\sigma \equiv 1) \log \frac{1}{\mu(\sigma \equiv 1)} \leq C_{LS}(\pi) \frac{\mu(\sigma \equiv 1)}{4}. 
$$
Hence,
$$
C_{LS}(\pi) \geq 4 \log \frac{1}{\mu(\sigma \equiv 1)} = 4 \log (2) n .
$$
For the complete graph, as for the previous examples, we get the right order of magnitude:
$$
4 \log (2) n \leq C_{LS}(\pi) \leq 6n .
$$

\subsection{The hypercube}
Consider the hypercube $G=\{0,1\}^N$. Assume that the lamplighter performs a simple random walk: each point has $N$ neighbors, $p(x,y)=1/N$ for $x \sim y$
($x \sim y$ if and only if 
$x$ and $y$ differ only in one point, or equivalently if and only if the Hamming distance between $x$ and $y$ is exactly 1), $\nu \equiv 1/2^N$. Assume also that $\mu$ is the uniform measure on $\{0,1\}^{\{0,1\}^N}$,
i.e. $\mu \equiv 1/2^{2^N}$, and that $c_x(\sigma)=1/2$.
Since the gap of the random walk on the hypercube is of order $N$, Theorem
\ref{th:main-intro} gives
\begin{equation} \label{eq:11}
C_{LS}((\{0,1\}^N)^\diamond ,\pi) \leq C N 2^N
\end{equation}
for some constant $C$ independent on $N$.

The hypothesis $(H)$ is not satisfied. As for the example of the complete graph,  one has to find a nice test function
in order to estimate from below the log-Sobolev constant.

Let $o:=(0,\ldots,0)$ be one "corner" of $G$, and define
$$
A_k:=B(o,\frac N2 - k) = \left\{ (x_1,\ldots,n_N) \in G : \sum_{i=1}^N x_i \leq \frac N2 - k
\right\}
$$
for $k=0,\ldots, N/2$ where for simplicity we will assume that $N$ is even.
Using e.g. the Stirling formula, one can easily see that there exists a constant
$C$ (independent on $N$) such that
$\left( \genfrac{}{}{0pt}{}{N}{N/2} \right) \leq C2^N / \sqrt N$.
Hence, Since $\left( \genfrac{}{}{0pt}{}{N}{i} \right)$ is increasing in $i$ for
$i \leq N/2$,
\begin{eqnarray}
|A_k| &=& \sum_{i=0}^{\frac N2 - k} \left( \genfrac{}{}{0pt}{}{N}{i} \right)
= \frac 12 2^N - \sum_{i=\frac N2 - k +1}^{N/2} \left( \genfrac{}{}{0pt}{}{N}{i} \right)
\nonumber \\
& \geq &
\frac 12 2^N - k \left( \genfrac{}{}{0pt}{}{N}{N/2} \right)
\geq \left( \frac 12 - \frac{kC}{\sqrt N} \right) 2^N \nonumber \\
& \geq &
\frac 14 2^N \label{eq:1}
\end{eqnarray}
for any $k=0, \ldots, \lfloor \delta  \sqrt N \rfloor$,
where $\delta := 1/(4C)$ is a parameter and $\lfloor \cdot \rfloor$ denotes the integer part.
Set $\delta_N:=\lfloor \delta  \sqrt N \rfloor$.

Our test function is
$$
f(\sigma,x) = \ind_{\sigma_{A_{\delta_N}} \equiv 1} (\sigma)
\sum_{k=0}^{\delta_N -2} \ind_{A_k^c}(x)
$$
where $\ind_{\sigma_{A_{\delta_N}} \equiv 1}$
is the indicator function of the event that all the lamps are on on the set
$A_{\delta_N}$, and $A_k^c$ denotes the complement of $A_k$.

Since $A_{\delta_N} \cap A_k^c = \emptyset$ for any
$k=0,\ldots,\delta_N -2$, the only non trivial contribution in the Dirichlet form
comes from the walk on $G$. Namely,
\begin{eqnarray*}
\mathcal{E}_\pi(f) & = &
\frac 14 \sum_\sigma \sum_x \sum_{y \sim x} \mu(\sigma) \nu(x) p(x,y)
(f(\sigma,y) - f(\sigma,x))^2 \\
& = &
\frac 1N \mu \left( \left\{ \sigma_{A_{\delta_N}} \equiv 1 \right\} \right)
\sum_x \sum_{y \sim x} \nu(x)
\left( \sum_{k=0}^{\delta_N -2} \ind_{A_k}(x) - \ind_{A_k}(y) \right)^2 \\
& \leq &
C \mu \left( \left\{ \sigma_{A_{\delta_N}} \equiv 1 \right\} \right) .
\end{eqnarray*}
for some constant $C$ depending on $\delta$.
In the second line we used the fact that $\ind_{A_k} = 1-\ind_{A_k^c}$, while the last
inequality comes from explicit counting.

On the other hand, note
that $\pi(f^2)=\mu \left( \left\{ \sigma_{A_{\delta_N}} \equiv 1 \right\} \right)
\nu(g^2)$
with
$g(x):=\sum_{k=0}^{\delta_N -2} \ind_{A_k^c}(x)$. Hence, a simple computation
and the fact that the entropy is a non negative function yield
\begin{eqnarray*}
\ent_\pi(f^2) & = &
\mu \left( \left\{ \sigma_{A_{\delta_N}} \equiv 1 \right\} \right)
\left( \nu(g^2) \log \frac{1}{\mu \left( \left\{ \sigma_{A_{\delta_N}} \equiv 1 \right\} \right)}
+ \ent_\nu(g^2) \right) \\
& \geq &
\nu(g^2) \mu \left( \left\{ \sigma_{A_{\delta_N}} \equiv 1 \right\} \right)
\log \frac{1}{\mu \left( \left\{ \sigma_{A_{\delta_N}} \equiv 1 \right\} \right)}
 .
\end{eqnarray*}
Since $g(x) = \delta_N$ on $A_0^c$, it is not difficult to see that there exists a constant
$C$ (independent on $N$) such that
$$
\nu(g^2) \geq \sum_{x \in A_{\delta_N}^c} \nu(x) g(x)^2
\geq C N .
$$
The very definition of the logarithmic Sobolev constant leads to
$$
C_{LS}((\{0,1\}^N)^\diamond ,\pi) \geq
C' N \log \frac{1}{\mu \left( \left\{ \sigma_{A_{\delta_N}} \equiv 1 \right\} \right)}
\geq C'' N 2^N
$$
since by \eqref{eq:1}, one has
$$
\mu \left( \left\{ \sigma_{A_{\delta_N}} \equiv 1 \right\} \right)
= \left( \frac{1}{2} \right)^{|A_{\delta_N}|}
\leq \left( \frac{1}{2} \right)^{2^N/4}
$$
where $C'$ and $C''$ are  constants independent on $N$.

By \eqref{eq:11}, we conclude that
$$
C_{LS}(\pi) = O( N 2^N ) .
$$
Peres-Revelle \cite{peres-revelle} result show that the uniform mixing time of the lamplighter on the hypercube is between $O(N2^N)$ and $O(N^2 2^N)$. Our estimates on the log-Sobolev constant give exactly the same bounds (and unfortunately no more).

\subsection*{Acknowledgements}
We warmly thank Yuval Peres and Amir Dembo for usefull discussions on the topic of this work and University
of Berkeley and Standford for their kind hospitality.


\bibliographystyle{alpha}
\bibliography{logsob-Lamplighter}

 \end{document}